\documentclass{article}
\usepackage{graphicx} 
\usepackage[utf8]{inputenc}
\usepackage{indentfirst}
\usepackage{amsmath}
\usepackage{amssymb}
\usepackage{amsthm}
\usepackage{esint}
\usepackage{graphicx}%
\usepackage{multirow}%
\usepackage{multicol}
\usepackage{amsmath,amssymb,amsfonts}%
\usepackage{amsthm}%
\usepackage{mathrsfs}%
\usepackage[title]{appendix}%
\usepackage{xcolor}%
\usepackage{textcomp}%
\usepackage{manyfoot}%
\usepackage{subcaption}%
\usepackage{booktabs}%
\usepackage{algorithm}%
\usepackage{algorithmicx}%
\usepackage{algpseudocode}%
\usepackage{listings}%
\usepackage{tikz}
\usetikzlibrary{decorations.markings}
\tikzstyle{vertex}=[circle, draw, inner sep=0pt, minimum size=6pt]

\usepackage[a4paper]{geometry}
\usepackage{optidef}

\theoremstyle{plain}
\newtheorem{thm}{Theorem}[section]

\newtheorem*{thm*}{Theorem}
\newtheorem*{lemma*}{Lemma}
\newtheorem*{prop*}{Lemma}
\newtheorem*{cor*}{Corollary}
\newtheorem*{conj*}{Conjecture}

\theoremstyle{definition}
\newtheorem{defn}[thm]{Definition}

\newtheorem*{defn*}{Definition}

\newtheorem{ques}[thm]{Question}

\theoremstyle{remark}

\DeclareMathOperator{\rank}{rank}

\usepackage{authblk}
\title{Maximum likelihood thresholds of Gaussian graphical models and graphical lasso}
\author{Daniel Irving Bernstein}
\author{Hayden Outlaw}
\affil{Tulane University}
\date{\today}

\begin{document}

\maketitle

\abstract{Associated to each graph $G$ is a Gaussian graphical model. Such models are often used in high-dimensional settings, i.e. where there are relatively few data points compared to the number of variables. The maximum likelihood threshold of a graph is the minimum number of data points required to fit the corresponding graphical model using maximum likelihood estimation. Graphical lasso is a method for selecting and fitting a graphical model.
In this project, we ask: when graphical lasso is used to select and fit a graphical model on $n$ data points, how likely is it that $n$ is greater than or equal to the maximum likelihood threshold of the corresponding graph?
Our results are a series of computational experiments.
}

\section{Introduction}

The set of all $p\times p$ positive definite matrices will be denoted $\mathcal{S}^p_{++}$.
Given a vector $\mu \in \mathbb{R}^p$ and positive definite matrix $K \in \mathcal{S}^p_{++}$,
the \emph{multivariate normal distribution} with mean vector $\mu$ and concentration matrix (i.e.~inverse covariance matrix) $K$ is the probability distribution with the following density function
\[
    f_{\mu,K}(x) := \frac{\exp(-\frac{1}{2}(x-\mu)^TK(x-\mu))\sqrt{\det(K)}}{(2\pi)^{p/2}}.
\]
A \emph{Gaussian model} is a set of Gaussian distributions.
Since $f_{\mu,K}(x) = f_{0,K}(x-\mu)$,
we will only concern ourselves with Gaussian distributions with densities of the form $f_{0,K}$.
Thus we can identify each Gaussian model with the set of allowable concentration matrices.
That is, a Gaussian model is simply a set of positive definite matrices.

Given a Gaussian model $\mathcal{M}\subseteq \mathcal{S}^p_{++}$ and a dataset $X \in \mathbb{R}^{p\times n}$
whose columns $x_1,\dots,x_n$ are believed to be independently and identically distributed (i.i.d.) according to a distribution in $\mathcal{M}$,
the corresponding \emph{maximum likelihood estimator (MLE)} is the solution to the following optimization problem
\begin{maxi}|l|
    {K}{\prod_{i=1}^n f_{0,K}(x_i)}{}{\label{eq: mlt optimization from distribution}}
    \addConstraint{K \in \mathcal{M}}.
\end{maxi}
The MLE is a statistically consistent estimator of the density function~\cite{casella2021statistical} but it does not always exist.
When it does exist, it is the solution to the following optimization problem,
which is convex when $\mathcal{M}$ is convex.
One can see this by taking the logarithm of the objective function in~\eqref{eq: mlt optimization from distribution} and doing some algebraic manipulation.
\begin{mini}|l|
   {K}{{\rm Tr}\left(XX^TK\right) - \log \det K}{}{\label{eq: mlt optimization}}
   \addConstraint{K \in \mathcal{M}}.
\end{mini}
For a given Gaussian model $\mathcal{M} \subseteq \mathcal{S}^p_{++}$,
the minimum $n$ such that the MLE exists for almost every $X \in \mathbb{R}^{p\times n}$ is called the \emph{maximum likelihood threshold (MLT)} of $\mathcal{M}$.

Given a graph $G$ on vertex set $1,\dots,p$ and edge set $E$,
the corresponding \emph{Gaussian graphical model} is the Gaussian model
\[
    \mathcal{M}_G := \{K \in \mathcal{S}^n_{++} : K_{ij} = 0 \ {\rm whenever} \ ij \notin E\}.
\]
We will abuse terminology and refer to the maximum likelihood threshold of $\mathcal{M}_G$ as the maximum likelihood threshold of $G$.
Gaussian graphical models were introduced by Dempster in~\cite{dempster1972covariance} for the purpose of fitting a Gaussian to data in the high-dimensional setting, i.e.~when there are more random variables than datapoints.
In light of this, determining the MLT of a graph is an important problem.

An upper bound on the MLT of a graph in terms of a related matrix completion problem can be derived from Dempster's original work (see e.g.~\cite{gross2018maximum})
which immediately implies that the MLT of a graph is at most the number of vertices.
Buhl showed that the MLT of a graph lies between its clique number and the clique number of a minimal chordal cover (i.e.~the treewidth plus one)~\cite{buhl1993existence}.
Uhler gave an easily computable upper bound for the MLT of a graph~\cite{uhler2011geometry}, shown not to be sharp by Blekherman and Sinn~\cite{blekherman2019maximum}.
Gross and Sullivant showed how Uhler's bound on the MLT of a graph $G$ can be understood in terms of rigidity-theoretic properties of $G$
and used this to determine the MLTs of some families of graphs.
Bernstein, Dewar, Gortler, Nixon, Sitharam and Theran built upon work of Gross and Sullivant to show how the MLT of a graph $G$ can be defined purely in terms of rigidity theoretic properties of $G$
and used this to determine the MLTs of several more families of graphs~\cite{bernstein2021maximum}
and further built upon this to determine the MLTs of all graphs with nine or fewer vertices~\cite{bernstein2022computing}.

Graphical lasso~\cite{friedman2008sparse} is a method that computes and fits a graphical model to a given dataset.
Given a dataset $X \in \mathbb{R}^{p\times n}$ and a regularization parameter $\alpha > 0$,
the \emph{graph lasso estimator} which we denote $K_\alpha(X)$, is the solution to the following optimization problem.
\begin{mini}|l|
   {K}{{\rm Tr}\left(XX^TK\right) - \log \det K + \alpha \sum_{i \neq j} |K_{ij}|}{}{\label{eq: graphical lasso optimization}}
   \addConstraint{K \in \mathcal{S}^n_{++}}.
\end{mini}
The graph lasso estimator can be interpreted as the maximum \emph{a posteriori} estimator of a certain model~\cite{wang2012bayesian}.
It is known to be statically inconsistent in some settings~\cite{heinavaara2016inconsistency}, but it always exists~\cite{ravikumar2011high-dimensional}.

Note that~\eqref{eq: graphical lasso optimization} is obtained from~\eqref{eq: mlt optimization} by setting $\mathcal{M} = \mathcal{S}^n_{++}$ and adding a penalty term for non-zero off-diagonal entries.
This is done to force $K_\alpha(X)$ to become sparser as $\alpha$ increases.
We let $G_{\alpha}(X)$ denote the graph on vertex set $\{1,\dots,p\}$ where $ij$ is an edge whenever $K_{\alpha}(X)_{ij} \neq 0$.
Then $K_{\alpha}(X) \in \mathcal{M}_{G_\alpha(X)}$ so one can view graph lasso as an algorithm that both selects and fits a Gaussian graphical model to a dataset.

Given a dataset $X \in \mathbb{R}^{p\times n}$, we aim to understand how likely it is that the MLT of $G_\alpha(X)$ is less than or equal to $n$.
In other words, how often does graph lasso select a model whose MLT is at most the number of data points used to find it?
More precisely, we ask the following.

\begin{ques}\label{ques:main question}
    Let $X \in \mathbb{R}^{p \times n}$ have columns that are i.i.d.~according to the Gaussian distribution with mean zero and identity covariance.
    Let $q_{p,n,\alpha}$ denote the probability that the MLT of $G_\alpha(X)$ is at most $n$.
    How does $q_{p,n,\alpha}$ depend on $p$, $n$, and $\alpha$?
\end{ques}

In this paper, we aim to answer Question~\ref{ques:main question} empirically for each pair $(p,n)$ with $p = 3,\dots,9$ and $n = 1,\dots,p$.








\section{Computing the MLT of a small graph}
The MLT of a graph $G$ on $p$ vertices is the minimum $n$ such that for every generic $p\times p$ positive semidefinite matrix $S$ of rank $n$,
there exists $T \in \mathcal{S}^p_{++}$ such that $T_{ij} = S_{ij}$ for all edges $ij$ of $G$~\cite{gross2018maximum}.
This can be phrased in the first order logic over the reals, so in principle, the MLT of a graph can be computed via a quantifier elimination algorithm~\cite{caviness2012quantifier}.
Such algorithms are prohibitively slow for this purpose.
For graphs on $p\le 9$ vertices, rigidity theory techniques can be used to compute the MLT exactly~\cite[Theorem 1.5]{bernstein2022computing};
this is given in Algorithm~\ref{alg:MLT}.
To describe this, we need some definitions.
Given a positive integer $m$, we use the shorthand $[m]$ to denote $\{1,\dots,m\}$.

\begin{defn}\label{defn:rigidity matrix}
    Let $G$ be a graph on vertex set $[p]$ and edge set $E$.
    Let $n$ a positive integer and let $X \in \mathbb{R}^{p\times n}$.
    The corresponding \emph{rigidity matrix} $R_n^G(X)$ is the matrix with rows indexed by $E$ and columns indexed by $[p]\times [n]$ given as follows
    \[
        R_n^G(X)_{e,(i,k)} =
        \begin{cases} 
            X_{ik} - X_{jk} & {\rm if} \ i,j \ \textnormal{are the vertices incident to } e \\
            0 & {\rm otherwise.}
        \end{cases}
    \]
\end{defn}

Assuming that $X \in \mathbb{R}^{p\times n}$ is \emph{generic}, the rank of $R_n^G(X)$
is independent of $X$. More precisely, there exists an integer $r_n^G$ such that if $X$ is sampled from $\mathbb{R}^{p\times n}$ from a probability distribution whose density function is mutually absolutely continuous with respect to Lebesgue measure (e.g.~a normal distribution),
then $\rank(R_n^G(X)) = r_n^G$ almost surely.
The minimum $n$ such that $r_{n-1}^G = |E|$ is called the \emph{generic completion rank (GCR)} of $G$.

\begin{thm}[{\cite{bernstein2022computing,gross2018maximum,uhler2011geometry}}]\label{thm:gcr mlt}
    Let $G$ be a graph.
    Then the MLT of $G$ is at most the GCR of $G$.
    When $G$ has 9 or fewer vertices, these quantities are equal.
\end{thm}

The inequality from Theorem~\ref{thm:gcr mlt} was proven by Uhler in~\cite{uhler2011geometry}.
Gross and Sullivant described its rigidity-theoretic interpretation in~\cite{gross2018maximum}.
That these quantities are equal for nine or fewer vertices was shown in~\cite{bernstein2022computing}.
The inequality can be strict for ten or more vertices~\cite{blekherman2019maximum},
but numerical simulations suggest that the GCR and MLT of an Erd\"os-Reyni random graph are equal with high probability~\cite{bernstein2021maximum}.

For our purposes, the upshot of Theorem~\ref{thm:gcr mlt} is that there exists an efficient randomized algorithm that computes the MLT of a graph on nine or fewer vertices. In particular, if $X\in \mathbb{R}^{p\times n}$ is chosen uniformly at random from the set of matrices with floating-point entries in $[-1,1]$,
then $\rank(R_n^G(X)) = r_n^G$ with high probability.
This gives us a procedure for computing the GCR of any graph that is correct with high probability -- see Algorithm~\ref{alg:MLT}.

\begin{algorithm}
\caption{Given an positive integer $p \le 9$ and a graph $G$ on vertex set $[p]$ and edge set $E$, compute the MLT of $G$}
\label{alg:MLT}
\begin{algorithmic}
\If {$E = \emptyset$}
\Return{1}
\Else
    \State $n \gets 1$ \Comment{Counter that increases until MLT is reached}
    \State $r \gets 0$ \Comment{Keeps track of the rank of the rigidity matrix}
    \While {$r < |E|$}
    \State Construct $X \in \mathbb{R}^{p\times n}$ by sampling each entry from the uniform distribution on $[-1,1]$
    \State $r \gets \rank(R_n^G(X))$
    \EndWhile
    
    \Return{$n$}
\EndIf
\end{algorithmic}
\end{algorithm}

\section{Experiments and results}
\begin{figure}[h]
\begin{multicols}{3}
\includegraphics[width = \linewidth]{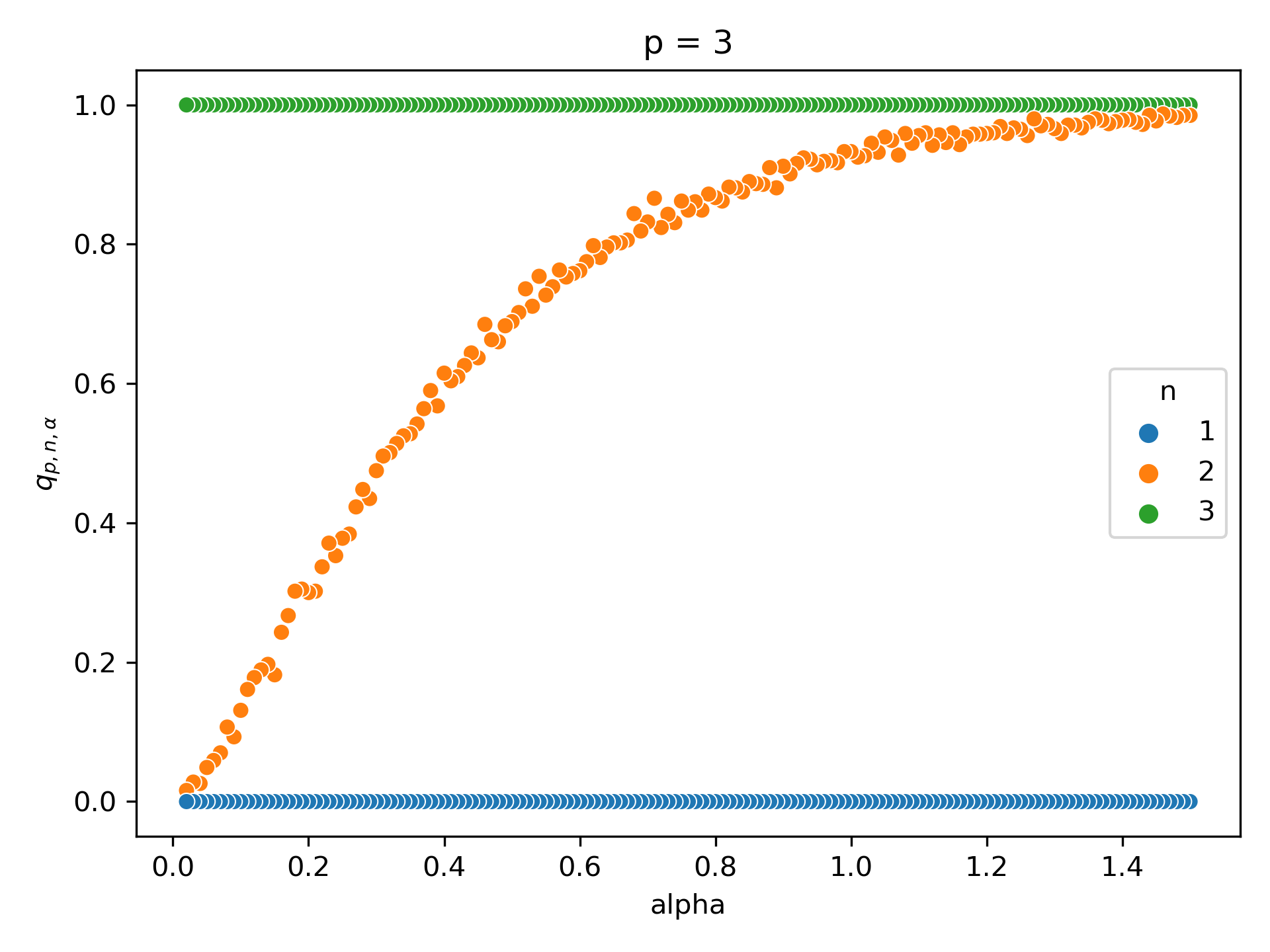}
\includegraphics[width = \linewidth]{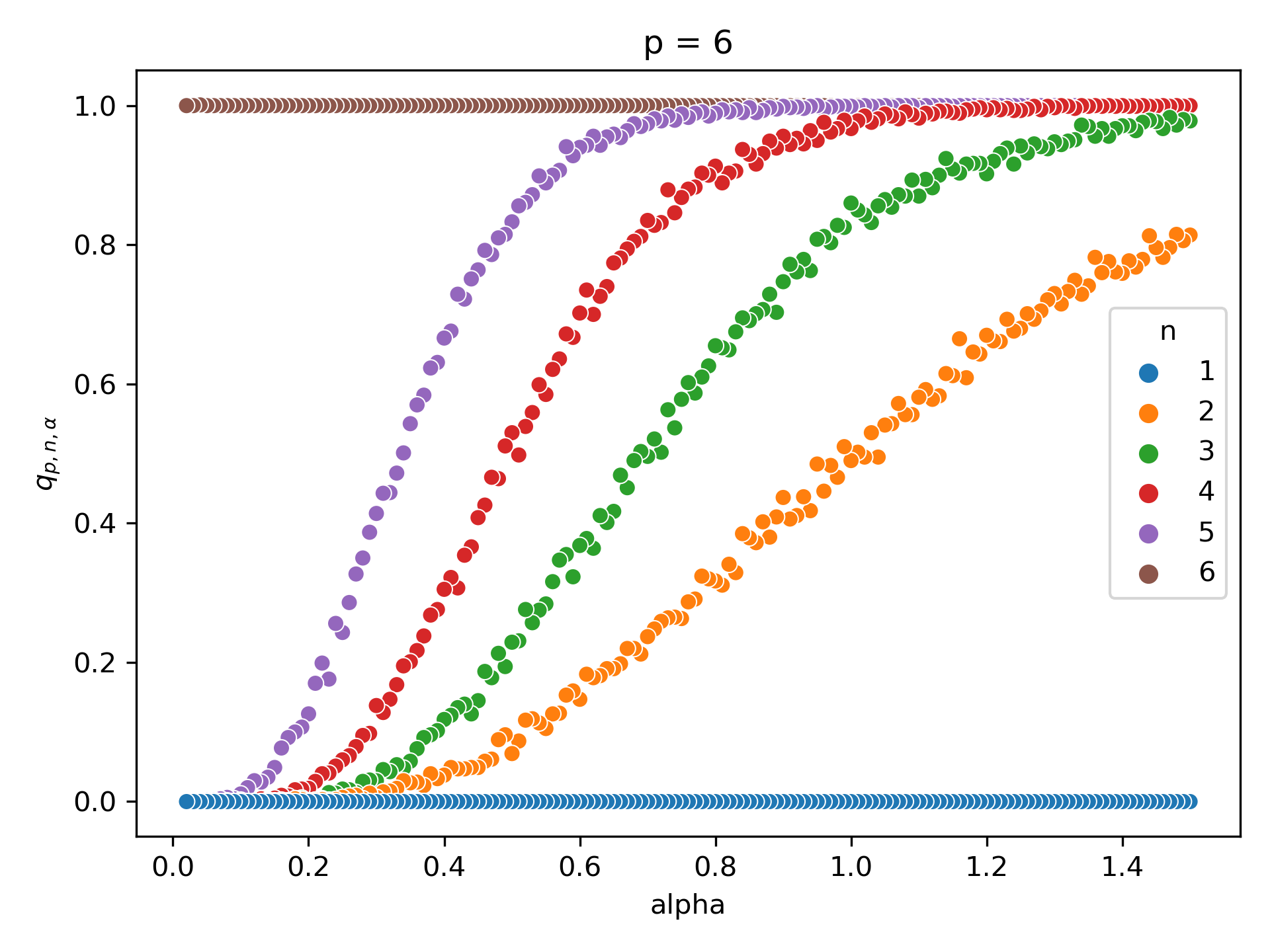}
\includegraphics[width = \linewidth]{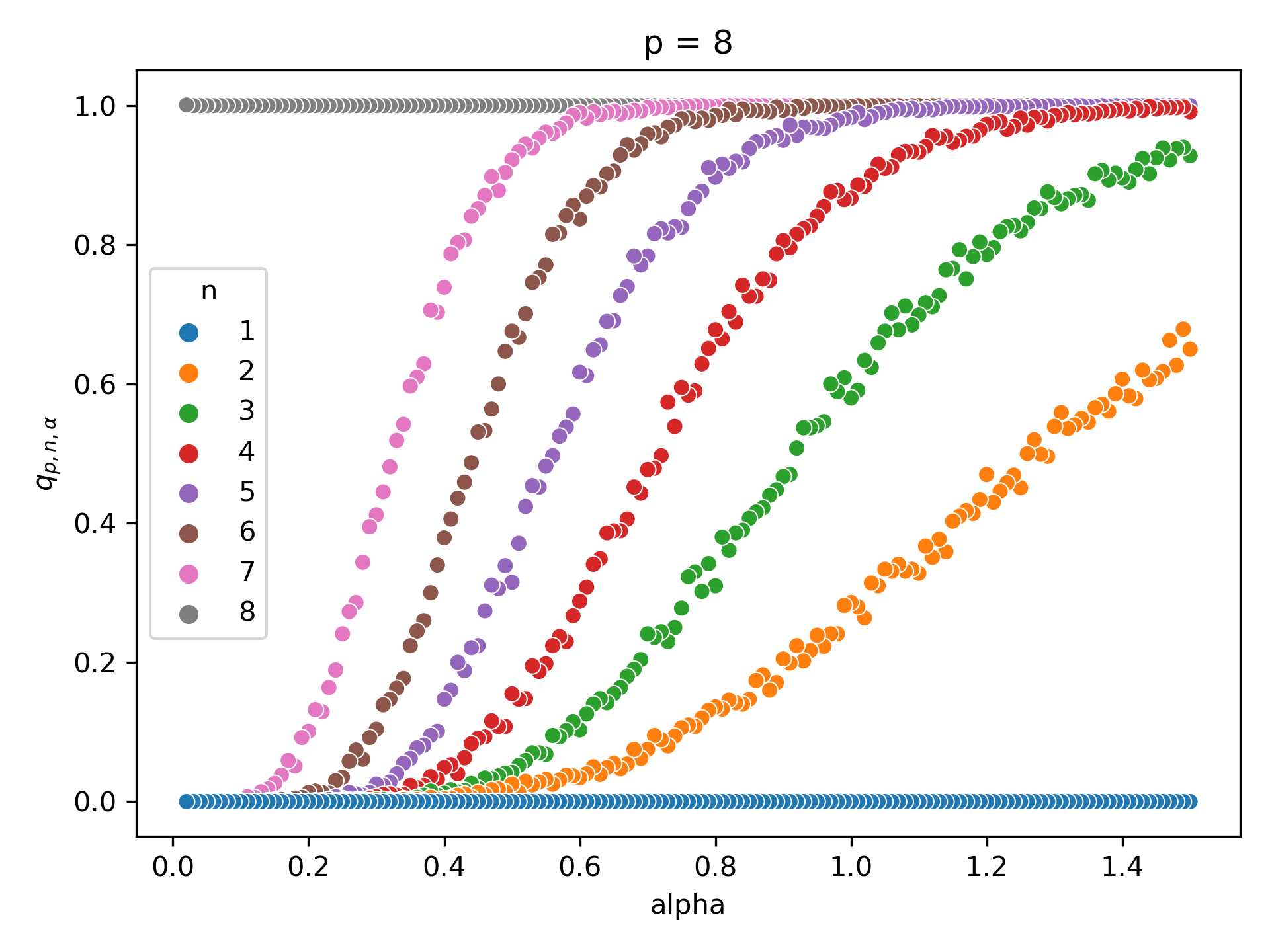}
\includegraphics[width = \linewidth]{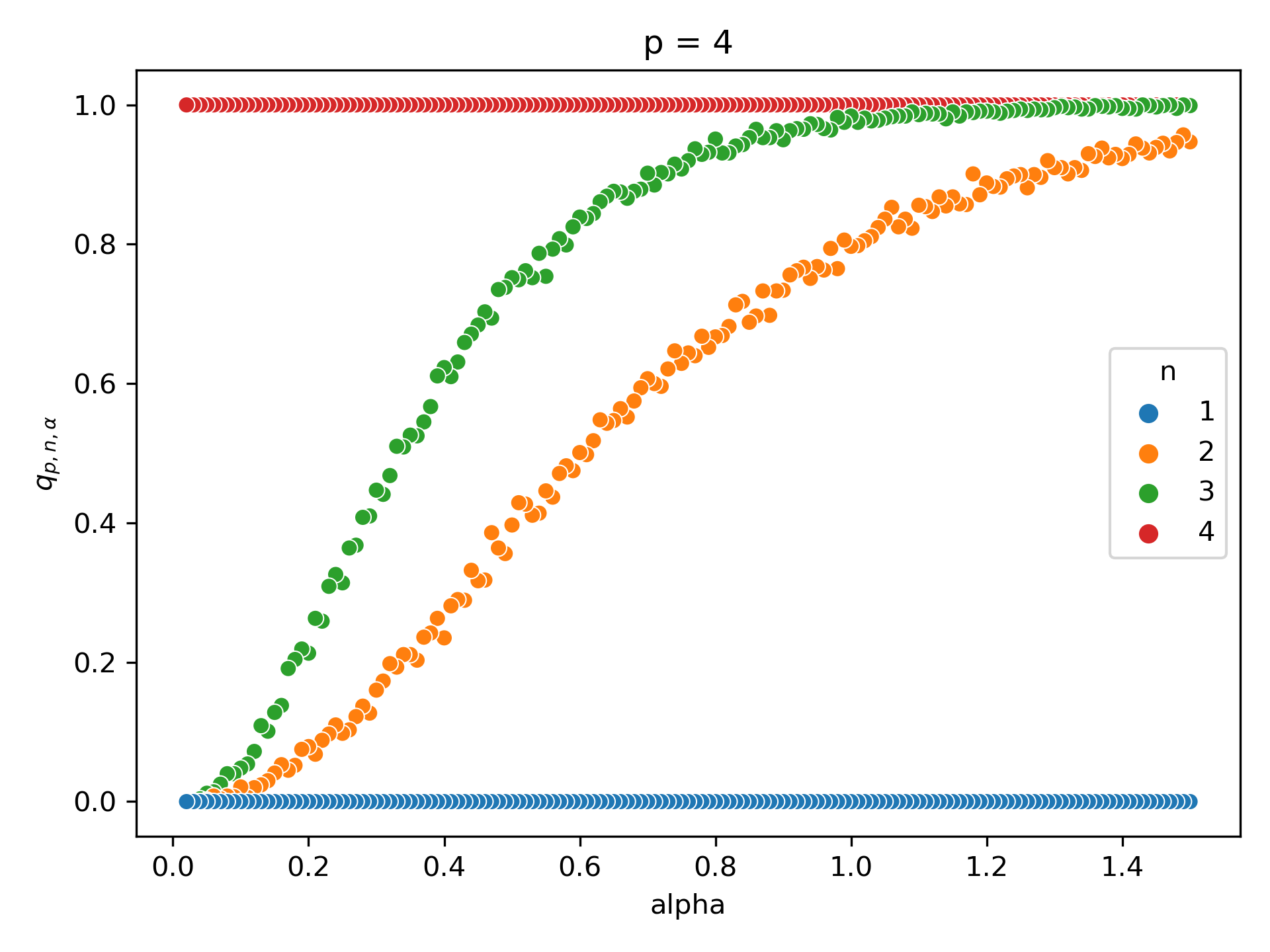}
\includegraphics[width = \linewidth]{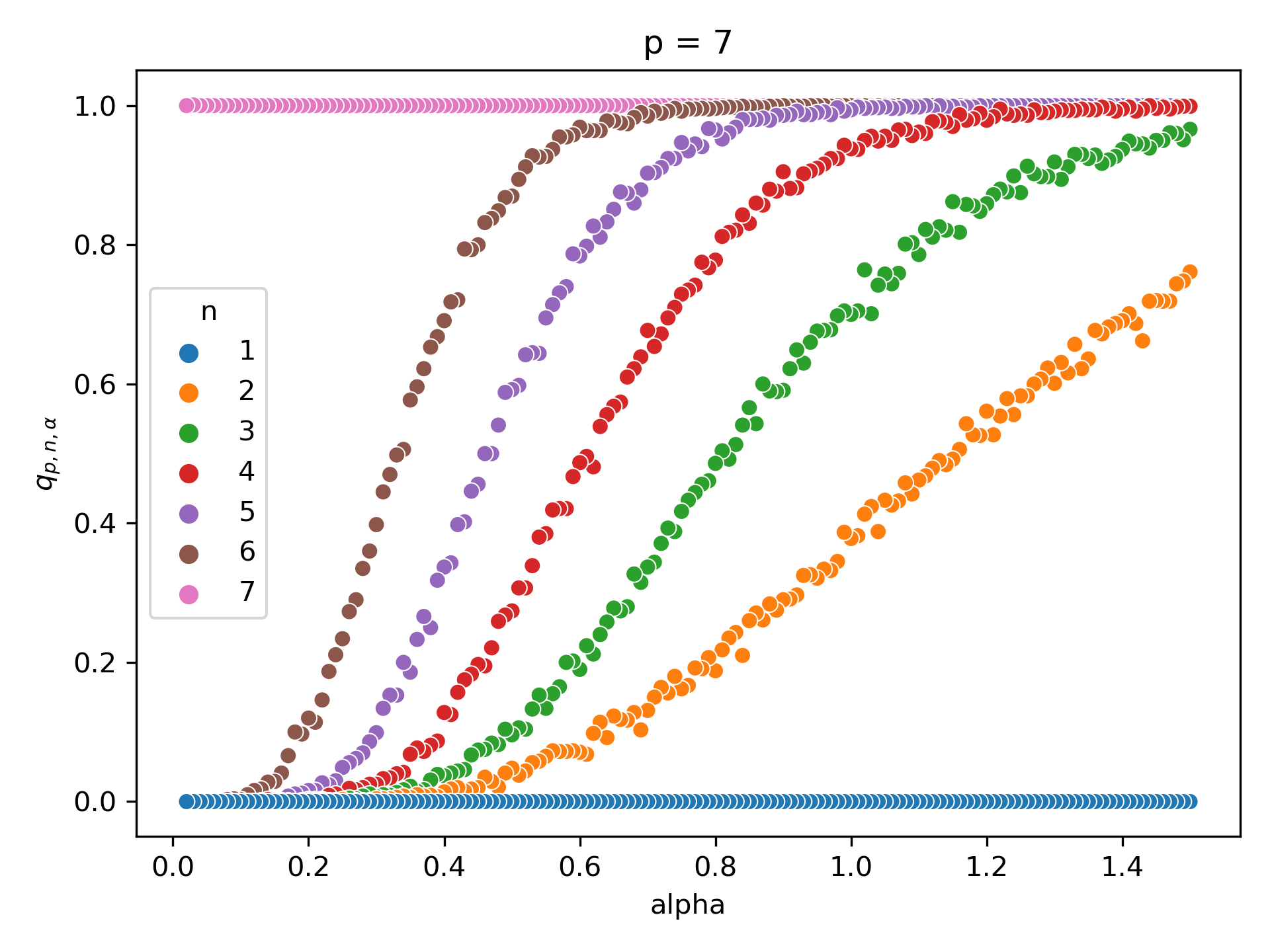}
\includegraphics[width = \linewidth]{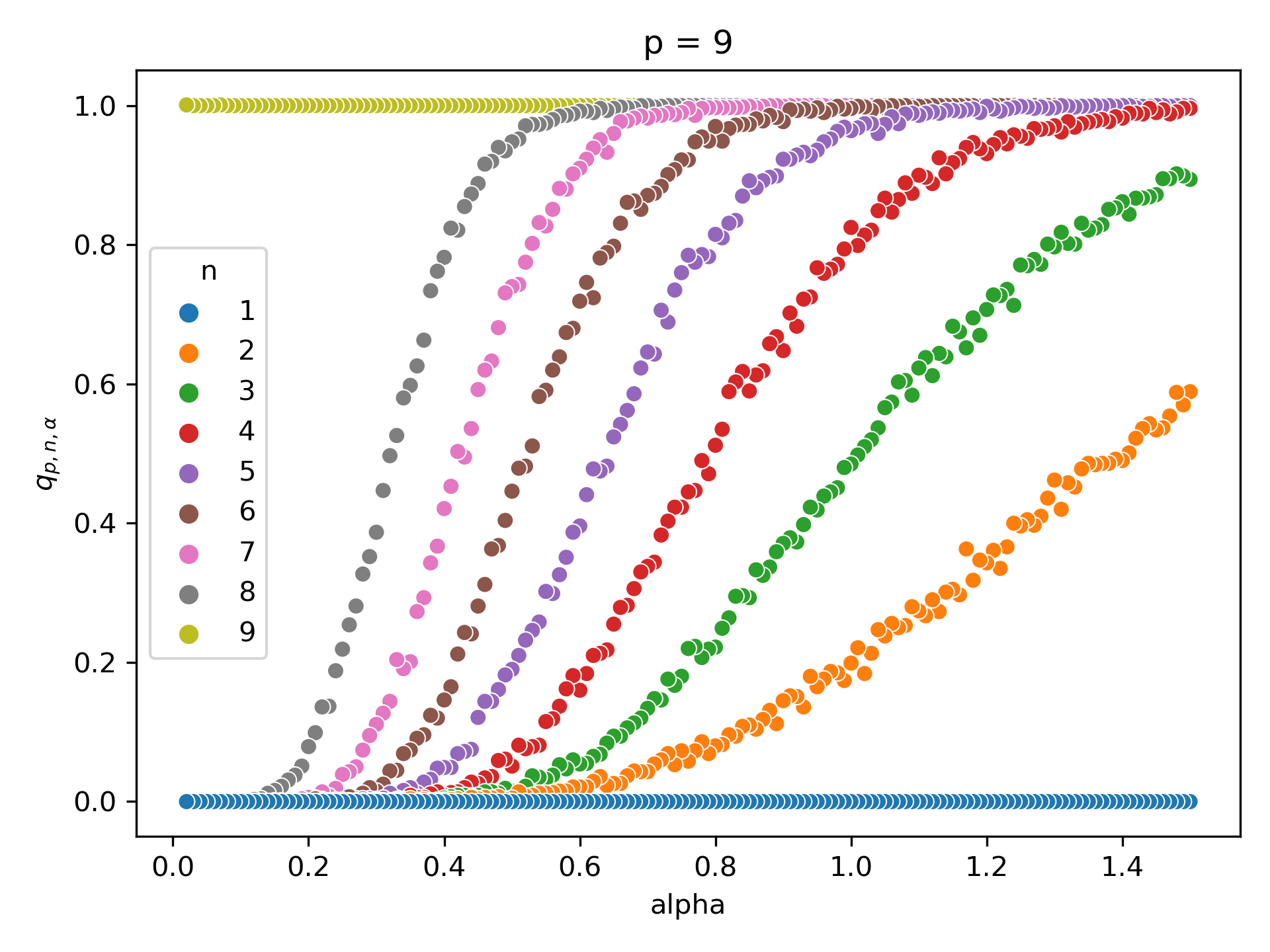}
\includegraphics[width = \linewidth]{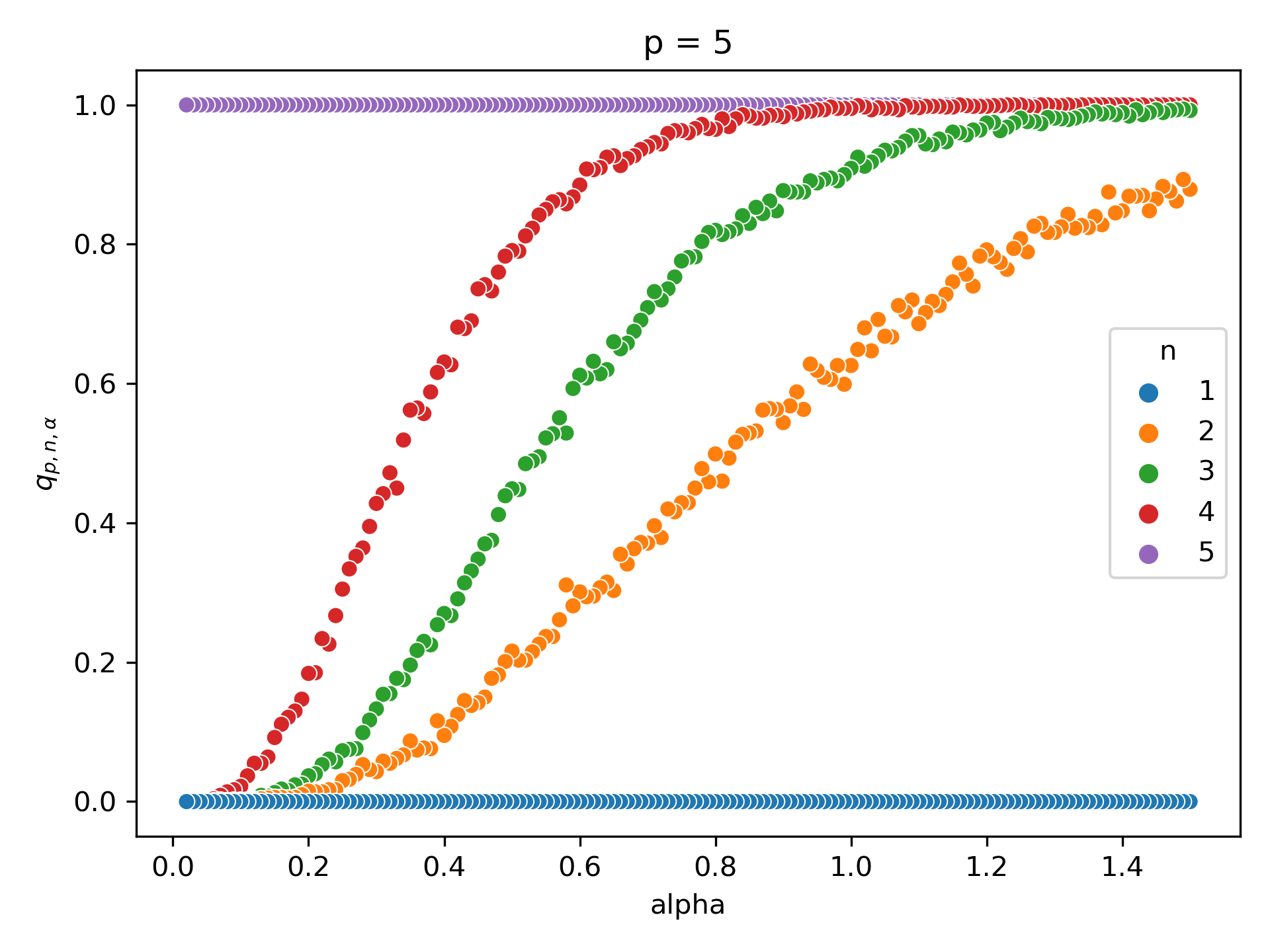}
\end{multicols}
\caption{Numerical estimation of the probability $q_{p,n,\alpha}$ that the MLT of the graphical model selected by graphical lasso on $p$ variables with regularization parameter $\alpha$ has MLT bounded above by $n$, the number of samples used.}
\label{fig:experimental results}
\end{figure}

We recall some definitions.
Let $K_\alpha(X)$ denote the output of the graphical lasso algorithm when applied to
the columns of $X$ with regularization parameter $\alpha$,
and let $G_\alpha(X)$ denote the graph that has $ij$ as an edge if and only if $K_\alpha(X)_{ij} \neq 0$.
Let $q_{p,n,\alpha}$ denote the probability that the MLT of $G_{\alpha}(X)$ is $n$ or less when
$X$ is a $p\times n$ matrix whose columns are i.i.d.~from the normal distribution on $\mathbb{R}^p$ with mean zero and identity covariance.
Figure~\ref{fig:experimental results} shows our experimental results that model how $q_{p,n,\alpha}$
depends on $p,n,\alpha$ for $3 \le p \le 9$ and $1 \le n \le p$.

For each fixed value of $(p,n)$, we let $\alpha$ range from $0.01$ to $1.5$ in steps of size $0.01$.
For each fixed value of $(p,n,\alpha)$ considered,
we sampled $X \in \mathbb{R}^{p\times n}$ $1000$ times and computed $G_\alpha(X)$ for each using \verb|sklearn.covariance.GraphicalLasso|~\cite{scikit-learn}, an implementation of graphical lasso in Python.

Let $\hat q_{p,n,\alpha}$ denote the proportion of times that the MLT of $G_{\alpha}(X)$ was less than $n$.
Since this is the proportion of successes obtained from $1000$ samples of a Bernoulli random variable with success probability $q_{p,n,\alpha}$,
$\hat q_{p,n,\alpha}$ is a statistically consistent estimator of $q_{p,n,\alpha}$ and the endpoints of the
$95 \% $ confidence interval are $\hat{p} \pm 1.96\sqrt{\frac{\hat{p}(1-\hat{p})}{1000}}$.
The results of these experiments are shown in Figure~\ref{fig:experimental results}.
No error bars are displayed since their range is too small to be visible.

\section{Discussion}
The MLT of a graph on $p$ vertices is never greater than $p$, which is reflected by all our numerical simulations.
It is interesting that graphical lasso never selected the graphical model with no edges,the only model with MLT $1$, when only one sample was used.

For each positive integer $k$, the $k$-core of a graph is the induced subgraph one obtains by repeatedly deleting vertices of degree strictly less than $k$.
Gross and Sullivant observed that the smallest $k$ for which the $k$-core of $G$ is empty is an upper bound on the MLT of $G$~\cite[Theorem 3.7]{gross2018maximum}.
As $k$ increases and the number of edges in $G$ decreases, the size of the $k$-core of $G$ decreases.
For this reason, one would expect sparse graphs to have empty $k$-cores for large values of $k$.
As the regularization parameter $\alpha$ increases,
the number of edges in $G_{\alpha}(X)$ decreases,
so one would expect $q_{p,n,\alpha}$ to increase as $\alpha$ increases, and this is reflected in our results.

For fixed values of $p$ and $\alpha$, as $n$ increases, $q_{p,n,\alpha}$ also increases.
On the one hand, this is maybe not surprising as $n$ increasing would make it easier for the MLT of $G$ to be less than $n$.
However, this is only if the expected MLT of $G_{\alpha}(X)$ increases by less than $1$ when $n$ increases by $1$.
It would be interesting to understand this phenomenon.

\bibliographystyle{plain}
\footnotesize
\bibliography{graphlasso}

\end{document}